\documentclass[pre,aps]{revtex4}
\usepackage{graphicx,epsf}
\usepackage{amssymb}
\usepackage{amsmath}
\usepackage{latexsym}

\def\be{\begin{equation}}
\def\ee{\end{equation}}
\def\bq{\begin{eqnarray}}
\def\eq{\end{eqnarray}}
\def\beq{\begin{eqnarray*}}
\def\eeq{\end{eqnarray*}}

\begin{document}

\title{Chaos in a predator-prey-based mathematical model for illicit drug consumption}

\author{Jean-Marc Ginoux$^{1,2}$, Roomila Naeck$^3$, Yusra Bibi Ruhomally$^4$, Muhammad Zaid Dauhoo$^4$, Matja{\v z} Perc$^{5,6,7}$}

\affiliation{$^1$Laboratoire LIS, CNRS, UMR 7020, Universit\'{e} de Toulon, BP 20132, F-83957 La Garde cedex, France}

\affiliation{$^2$Departament de Matem\`{a}tiques, Universitat Aut\`{o}noma de Barcelona, 08193 Bellaterra, Barcelona, Catalonia, Spain}

\affiliation{$^3$PSASS, Ecoparc de Sologne, Domaine de Villemorant, Neug-sur Beuvron 41210, France}

\affiliation{$^4$Department of Mathematics, University of Mauritius, Reduit, Mauritius}

\affiliation{$^5$Faculty of Natural Sciences and Mathematics, University of Maribor, Koro{\v s}ka cesta 160, SI-2000 Maribor, Slovenia}

\affiliation{$^6$Center for Applied Mathematics and Theoretical Physics, University of Maribor, Mladinska 3, SI-2000 Maribor, Slovenia}

\affiliation{$^7$Complexity Science Hub Vienna, Josefst\"{a}dterstra{\ss}e 39, A-1080 Vienna, Austria}

\begin{abstract}
Recently, a mathematical model describing the illicit drug consumption in a population consisting of drug users and non-users has been proposed. The model describes the dynamics of non-users, experimental users, recreational users, and addict users within a population. The aim of this work is to propose a modified version of this model by analogy with the classical predator-prey models, in particular considering non-users as prey and users as predator. Hence, our model includes a stabilizing effect of the growth rate of the prey, and a destabilizing effect of the predator saturation. Functional responses of Verhulst and of Holling type II have been used for modeling these effects. To forecast the marijuana consumption in the states of Colorado and Washington, we used data from Hanley (2013) and a genetic algorithm to calibrate the parameters in our model. Assuming that the population of non-users increases in proportion with the demography, and following the seminal works of Sir Robert May (1976), we use the growth rate of non-users as the main bifurcation parameter. For the state of Colorado, the model first exhibits a limit cycle, which agrees quite accurately with the reported periodic data in Hanley (2013). By further increasing the growth rate of non-users, the population then enters into two chaotic regions, within which the evolution of the variables becomes unpredictable. For the state of Washington, the model also exhibits a periodic solution, which is again in good agreement with observed data. A chaotic region for Washington is likewise observed in the bifurcation diagram. Our research confirms that mathematical models can be a useful tool for better understanding illicit drug consumption, and for guiding policy-makers towards more effective policies to contain this epidemic.
\end{abstract}

\maketitle

\section{Introduction}

Mathematical modeling of illicit drug consumption is a very difficult and complex problem. To this aim predator-prey models have been used at the end of the nineties. Then, in 2013 a model called NERA has been built to describe the dynamics of the non (N), experimental (E), recreational (R) and addict (A) user categories, respectively, within a given population. However, the original NERA model didn't involve limitation in drug consumption and was consequently unable to transcribe the periodic evolution of each category. So, we have modified this model by analogy with the classical predator-prey models and while considering non-users (N) as prey and users (E, R and A) as predator. Then, using data from the state of Colorado and Washington and a genetic algorithm, we calibrated our predator-prey NERA models by estimating their parameters sets. This allowed us to account for the periodic evolution of each category. Then, by considering that the population of Nonusers increases in proportion of the demography we highlighted chaotic regions within which the evolution of the variables becomes unpredictable. Thus, it appears that our validated NERA model can be a precious tool in forecasting of illicit drug consumption and can be of substantial interest to policy-makers in the problematic of illicit drug consumption.

Since the beginning of the nineties, many continuous time models have been proposed to describe the dynamics of drug consumption \cite{Baveja1, Baveja2, Caulkins1, Caulkins2, Caulkins3, Gragnani1, Gragnani3}. They mainly consisted in second order nonlinear models involving two variables which were used for predicting the long term dynamics of addicts and dealers of a large drug market, like that of an entire country. Few years later, Gragnani \textit{et al.} \cite{Gragnani2} proposed an extension to such models by adding a third ordinary differential equation to the original two-dimensional dynamical systems. This third variable represented the enforcement exerted by the authorities. Let's notice on the one hand that these models already used limitations in the growth and decay of each variable (at least for the two first) and on the other hand that, Gragnani \textit{et al.} \cite{Gragnani2} proved the existence of slow-fast limit cycles according to the singular perturbation theory as solution of their three-dimensional dynamical system.
Recently, Dauhoo \textit{et al.} \cite{Dauhoo} have considered that drug users are generally classified into three main categories, depending on their consumption frequency and the control they have over the drug, i.e., Experimental (E), Recreational (R) and Addicts (A) users. By adding a fourth variable, i.e., the Non (N) users to these first three ones, they proposed the NERA model. Although this four-dimensional dynamical system takes into account the mutual influence that drug users (E, R and A) can have on non-users (N) and on each other it does not contain any limitation in the growth and decay of each variable. Thus, no oscillatory or chaotic regime could be observed. That's the reason why we have proposed to modify Dauhoo's NERA model \cite{Dauhoo} by introducing a limitation in each ``functional response''. In their paper, Dauhoo \textit{et al.} \cite{Dauhoo} wrote the following sentence: ``Anyone could be a `prey' to illicit drugs''. This led us to make the analogy with predator-prey models.

\section{General Predator-Prey NERA Model}

This deterministic model aims at transcribing into mathematical functions variations of the number of individuals of each group due to their interactions with others. We make the assumption that such interactions are mainly characterized by the influence that individuals of one group may exert on the others. Thus, we consider that people influenced by a group leaves the group to which they belong to join the group which has influenced them. As a consequence, some individuals disappear of one group and appear in another. So, by analogy with the predator-prey models used for a long time in Theoretical Ecology \cite{Maybook,Scudo}, these influences, which cause such appearances and disappearances, \textit{i.e.}, increases and decreases (variations) of the number of individuals of each group can be regarded as predation of on group on another. In the predator-prey model that we propose, Non-users (N) represent prey for all other groups (E, R and A). Then,  Experimental-users (E) are predator of both Recreational-users (R) and Non-users (N) while Addict-users (A) are predator of all R, E and N-users. According to Kuznetsov \cite{Kuznetsov}, for the model to be realistic it is necessary to introduce a ``stabilizing effect of the competition among prey and a destabilizing effect of the predator saturation''. To this aim, we have used two different types of ``functional responses'' for the growth of prey (N) and for the growth of the predators (E, R, A). We have first considered that, in the absence of any predator, prey growth (N) must be limited. Such limitation or stabilizing effect is generally introduced by using the \textit{logistic law} introduced by Verhulst \cite{Verhulst}. Concerning the predators (E, R, A), the saturation of the predator rate (destabilizing effect) can be modeled with the classical Holling type II ``functional response'' \cite{Holling1,Holling2}. Of course, various other ``functional responses'' could have been chosen to this aim \cite{GinouxVolterra}. Let's notice that such saturation in the predator rate represents a limitation in the influence of each predator group on the others. At last, we consider that in the absence of its predators, the number of individuals of one group (E, R, and A except N) can decrease by a kind of ``natural mortality'' which can correspond to individuals leaving this group. In the most dramatic case of drug addict, this could be due to overdose.

\subsection{Functional responses}
\label{Sec1}

In 1837, the Belgian biologist Pierre-Fran\c{c}ois Verhulst proposed a model that took into account the limitation imposed by the increasing population size of a prey called X in absence of any predator. This model, called \textit{logistic law}, can be written as follows:

\[
\frac{dX}{dt} =  \beta X\left( 1 - X \right)
\]

In 1926, the Italian mathematician Vito Volterra developed the very first predator-prey models. The formulation of the equation representing the predation is based on the \textit{m\'{e}thode des rencontres} (method of encounters) and on the hypoth\`{e}se des \'{e}quivalents (hypothesis of equivalents) \cite{vol26, vol27, vol28, vol31}. The former assumes that for predation to occur between a predatory species (X) and a prey species (Y), it is first necessary to have encounters between these two species and that the number of encounters between them is proportional to the product of the number of individuals composing them, that is $\alpha X(t)Y(t)$, the coefficient of proportionality $\alpha$ being equal to the probability of an encounter. The second hypothesis consists in assuming that ``there is a constant ratio between the disappearances and appearances of individuals caused by the encounters'', that is, that predation of the preys is equivalent to increase of the predators. At the beginning, Volterra considers this increase as immediate. This means that predation is immediately transcribed in terms of growth of the predator species, whereas its effect naturally occurs with some delay. In his works Volterra \cite{vol31} also took this delay into account. This won't be the case in this paper. In the following, we will consider that the effects of influence that individuals of one group (X) may exert on another (Y) are analogous to the effects of predation of (Y) on (X). The mathematical function corresponding to the modeling of a behavior such as influence or predation is called a ``functional response''. The functional response proposed by Volterra to describe predation does not take into account any limitation, \textit{i.e.} ``satiety'' of the predator and so, the predation rate is a ``linear function'' of the prey. Thus, from the beginning of the thirties various types of ``function responses'' were proposed while using nonlinear mathematical functions presenting an asymptotic behavior and so a limitation \cite{gau,gau2}.

In the late 1950s, entomologist Crawford Stanley Holling \cite{Holling1, Holling2} developed two new functional responses for predation, also intended to
describe a certain satiety of the predator (Y) with respect to its prey (X): Holling function of type II and Holling function of type III. In this paper, we will only use the Holling type II functional response which can be represented by:

\[
\frac{X}{h + X}Y
\]

where $h$ represents half-saturation, that is, the value of the prey density $X = h$ for which the predation level reaches a value equal to half its maximum.

So, in this work we propose to used both \textit{logistic law} and Holling type II functional response for modeling the influence that exert the predators A, R and E on each others and on the prey N (See Fig. 1).

\begin{figure}[htbp]
\centerline{\includegraphics[width = 15.47cm,height = 20.94cm]{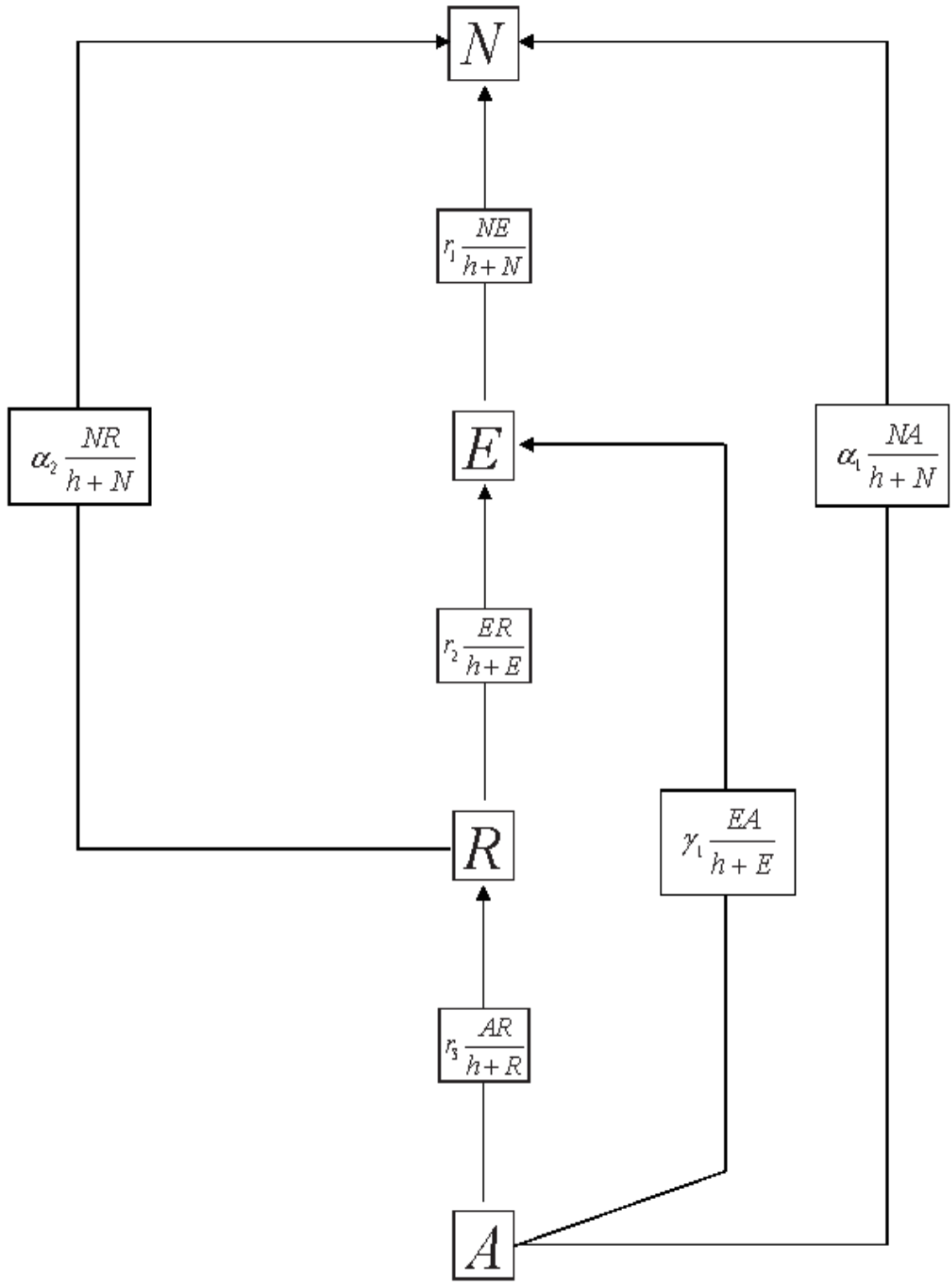}}
\caption{Schematic representation of the predator-prey NERA model.}
\label{fig1}
\end{figure}

\subsection{Model equations}

So, we have the following system of ordinary differential equations:

\begin{equation}
\label{eq1}
\begin{aligned}
\dfrac{dN}{dt} & = \beta_1 N \left( 1 - N \right) - r_1 \frac{N}{h+N}E - \alpha_1 \frac{N}{h+N}A - \alpha_2 \frac{N}{h+N}R, \hfill \\
\dfrac{dE}{dt} & = r_1 \frac{N}{h+N}E - r_2 \frac{E}{h+E}R - \beta_2 E - \gamma_1 \frac{E}{h+E}A, \hfill \\
\dfrac{dR}{dt} & = r_2 \frac{E}{h+E}R - \beta_3 R - r_3 \frac{R}{h+R}A + \alpha_2 \frac{N}{h+N}R, \\
\dfrac{dA}{dt} & = r_3 \frac{R}{h+R}A - \beta_4 A + \alpha_1 \frac{N}{h+N}A + \gamma_1 \frac{E}{h+E}A,
\end{aligned}
\end{equation}

where $\beta_1$ is the growth rate of the population of the prey ($N$) in the absence of any predator ($E,R,A$), $\beta_i$ with $i=2,3,4$ are the ``natural mortality'' of each predator ($E,R,A$) in the absence of all others and $r_i$ with $i=1,2,3$ is the predation rate of $A$ on $R$, $R$ on $E$ and $E$ on $N$ respectively. $\alpha_1$ and $\gamma_1$ represent the predation rate of $A$ on $N$ and $E$ respectively while $\alpha_2$ is that of $R$ on $N$. Thus, in this four-dimensional dynamical system, a set of 11 positive parameters: $(\beta_1, \beta_2, \beta_3, \beta_4, r_1, r_2, r_3, \alpha_1, \alpha_2, \gamma_1, h)$ is used.

\vspace{0.1in}

\textbf{Remark.} Let's notice that for each Holling type II functional response a different half saturation $h$ could have been chosen. Nevertheless, the aim of this work is to propose the most simple and consistent model for illicit drug consumption.

The sociological meaning of each parameter used in this model (\ref{eq1}) is given in Table \ref{tab1}.

\begin{table}[h]
	\centering
	\caption{Interpretation of the parameters in \textit{NERA} model}
	\begin{tabular}{ll}
		\hline\hline
		Parameter & Sociological Meaning \\ \hline
		$\displaystyle r_{1}$ & Influence rate of $\displaystyle E(t)$ on $\displaystyle N(t)$ \\
		
		$\displaystyle r_{2}$ & Influence rate of $\displaystyle R(t)$ on $\displaystyle E(t)$ \\
		$\displaystyle r_{3}$ & Rate at which recreational users change to addicts \\
		$\displaystyle \alpha_{1}$ & Influence rate of $\displaystyle A(t)$ on $\displaystyle N(t)$ \\
		$\displaystyle \alpha_{2}$ & Influence rate of $\displaystyle R(t)$ on $\displaystyle N(t)$ \\
		$\displaystyle \gamma_{1}$ & Influence rate of $\displaystyle A(t)$ on $\displaystyle E(t)$ \\
		$\displaystyle \beta_{1}$ & Rate of moving in and out of the Nonuser category \\
		$\displaystyle \beta_{2}$ & Rate at which experimental users quit drugs \\
			$\displaystyle \beta_{3}$ & Rate at which recreational users quit drugs \\
		$\displaystyle \beta_{4}$ & Rate at which addicts quit drugs \\
		\hline
	\end{tabular}
\label{tab1}
\end{table}	

\subsection{Dynamic aspects}

Due to the presence of the Holling type II functional responses, the determination of the fixed points of this four-dimensional dynamical system (\ref{eq1}) is not trivial while using the classical nullclines method. Nevertheless, by posing $E = R = A =0$ two obvious fixed points can be easily found:

\[
O(0, 0, 0,0) \quad \mbox{and} \quad I_1(1, 0, 0, 0).
\]

By posing $R = A =0$, a third fixed point can be also obtained:

\[
I_2 (N_2^*, E_2^*, 0, 0) = \left( \frac{\beta_2 h}{r_1 - \beta_2}, \frac{r_1 - \beta_2 \left( 1 + h \right)}{\left(r_1 - \beta_2 \right)^2} \beta_1 h  , 0, 0 \right)
\]

It follows that this fixed point $I_2$ is positive provided that:

\begin{equation}
\label{eq2}
r_1 - \beta_2 > 0 \mbox{ and } r_1 - \beta_2 \left( 1 + h \right) >0.
\end{equation}

Then, by posing $E = R =0$, a fourth fixed point can be also obtained:

\[
J_2 (N_2^*, 0, 0, A_2^*) = \left( \frac{\beta_4 h}{\alpha_1 - \beta_4},0,0, \frac{\alpha_1 - \beta_4 \left( 1 + h \right)}{\left(\alpha_1 + \beta_4 \right)^2} \beta_1 h \right)
\]

It follows that this fixed point $J_2$ is positive provided that:

\begin{equation}
\label{eq3}
\alpha_1 - \beta_4 > 0 \mbox{ and } \alpha_1 - \beta_4 \left( 1 + h \right) >0.
\end{equation}

Finally, while posing $A =0$, a fifth fixed point $I_3 (N_3^*, E_3^*, R_3^*, 0)$ can be determined but its expression is too long to be written here.

\vspace{0.1in}

\textbf{Remark.} Let's notice that all fixed points depend on parameters $(\beta_1)$.

\vspace{0.1in}

Following the works of Freedman \cite{Freed}, the system (\ref{eq1}) may be written as:

\begin{equation}
\label{eq4}
\begin{aligned}
\dfrac{dN}{dt} & = N g \left( N \right) - \left( r_1 E + \alpha_1 A + \alpha_2 R \right) p_1\left( N \right), \hfill \\
\dfrac{dE}{dt} & = E \left( - \beta_2 + r_1 p_1\left( N \right) \right) - \left(r_2 R + \gamma_1 A \right) p_2 \left( E \right), \hfill \\
\dfrac{dR}{dt} & = R \left(- \beta_3 + r_2 p_2 \left( E \right) + \alpha_2 p_1 \left( N \right) \right) - r_3 A p_3 \left( E \right) , \\
\dfrac{dA}{dt} & = A \left[ - \beta_4 + r_3 p_3 \left( R \right) + \alpha_1 p_1 \left( N \right) + \gamma_1 p_2 \left( E \right) \right],
\end{aligned}
\end{equation}

where $g(N) = \beta_1 (1 - N)$, $p_1(N) = \dfrac{N}{h + N}$, $p_2(E) = \dfrac{E}{h + E}$ and $p_3(R) = \dfrac{R}{h + R}$. Such a formulation will simplify the computation of the eigenvalues of the functional Jacobian matrix (\ref{eq5}) presented below.

\subsubsection{Functional Jacobian matrix}

The Jacobian matrix of system (\ref{eq4}) reads:

\begin{equation}
\label{eq5}
J = \begin{pmatrix}
J_{11} & -r_1 p_1 ( N ) &  - \alpha_2 p_1 ( N ) & - \alpha_1 p_1( N ) \vspace{6pt} \hfill \\
r_1 E p'_1 ( N ) & J_{22} &  - r_2 p_2 ( E ) & - \gamma_1 p_2( E ) \vspace{6pt} \hfill  \\
\alpha_2 R p'_1 ( N ) & r_2 R p'_2 ( E ) & J_{33} & - r_3 p_3( R ) \vspace{6pt} \hfill  \\
\alpha_1 A p'_1 ( N ) & \gamma_1 A p'_2 ( E ) &  - r_3 A p'_3 ( R ) & J_{44} \vspace{6pt} \hfill
\end{pmatrix}
\end{equation}

where the diagonal terms read:

\[
\begin{aligned}
J_{11}\left(N,E,R,A \right) & = g\left( N \right) + N g'\left( N \right) - \left( r_1 E + \alpha_1 A + \alpha_2 R\right) p'_1 \left( N \right), \\
J_{22}\left(N,E,R,A \right) & = - \beta_2 + r_1 p_1 \left( N \right) - \gamma_1 u p'_2 \left( E \right), \\
J_{33}\left(N,E,R,A \right) & = - \beta_3 + r_2 p_2 \left( E \right) + \alpha_2 p_1 \left( N \right) - r_3 A p'_3 \left( R \right), \\
J_{44}\left(N,E,R,A \right) & = - \beta_4 + r_3 p_3 \left( R \right) + \alpha_1 p_1 \left( N \right) + \gamma_1 p_2 \left( E \right).
\end{aligned}
\]

\vspace{0.1in}

Let's notice that for $i=1,2,3$ we have for $N=0$ and then for $N=1$:

\begin{equation}
\label{eq6}
\begin{aligned}
& g\left( 0 \right) = \beta_1 \mbox{, } g'\left( 0 \right) = -\beta_1 \mbox{, } p_i\left( 0 \right) = 0\mbox{, } p'_i\left( 0 \right) = \frac{1}{h},\\
&g\left( 1 \right) = 0\mbox{, }g'\left( 1 \right) = -\beta_1\mbox{, }p_i\left( 1 \right) = \frac{1}{h+1}\mbox{, }p'_i\left( 1 \right) = \frac{h}{\left(h+1\right)^2}.
\end{aligned}
\end{equation}

\subsubsection{Eigenvalues at O(0,0,0,0)}

Taking into account the above conditions (\ref{eq6}), the functional Jacobian matrix (\ref{eq5}) is diagonal. Thus, the four eigenvalues are:

\[
\lambda_1 = \beta_1 \mbox{ ; } \lambda_2 = - \beta_2 \mbox{ ; } \lambda_3 = - \beta_3 \mbox{ ; } \lambda_4 = - \beta_4.
\]

It follows that the origin $O$ is a \textit{saddle}.

\subsubsection{Eigenvalues at $I_{1}$(1,0,0,0)}

Taking into account the above conditions, the functional Jacobian matrix (\ref{eq5}) is diagonal. Thus, the four eigenvalues are:

\[
\lambda_1 = - \beta_1 \mbox{ ; } \lambda_2 = - \beta_2 + \frac{r_1}{h + 1} \mbox{ ; } \lambda_3 = - \beta_3 + \frac{\alpha_2}{h + 1} \mbox{ ; } \lambda_4 = - \beta_4 + \frac{\alpha_1}{h + 1}.
\]

According to conditions (\ref{eq2}-\ref{eq3}), both $\lambda_2$ and $\lambda_3$ are positive. So, whatever the values of the parameters, it follows that $I_1$ is also a \textit{saddle}.

\subsubsection{Eigenvalues at $I_{2}$($\rm{\textit{N}_2^*}$, $\rm{\textit{E}_2^*}$, 0, 0) and at $J_{2}$($\rm{\textit{N}_2^*}$, 0, 0, $\rm{\textit{A}_2^*}$)}

Although the four eigenvalues evaluated at $I_2$ and $J_2$ are too long to be expressed, two of them contain a square root. So, according to the choice of the parameters, these eigenvalues may be complex conjugate. Thus, if we assume that the expression in the square root is negative, we can look for the sign of the remaining part which can be considered as the \textit{real part} of these eigenvalues. Such a real part is positive provided that:

\[
0 < r_1 - \beta_2 \left( 1 + h \right) - h r_1 \quad \mbox{ and } \quad 0 < r_1 - \beta_2,
\]
\[
0 < \alpha_1 - \beta_4 \left( 1 + h \right) - h \alpha_1 \quad \mbox{ and } \quad 0 < \alpha_1 - \beta_4.
\]

Combining these conditions with the previous one (\ref{eq2}-\ref{eq3}), we find that

\begin{equation}
\label{eq7}
\begin{aligned}
& 0 < h r_1 < r_1 - \beta_2 \left( 1 + h \right),\\
& 0 < h \alpha_1 < \alpha_1 - \beta_4 \left( 1 + h \right).
\end{aligned}
\end{equation}

It follows that $I_2$ and $J_2$ are a \textit{saddle-foci} (two eigenvalues are complex conjugate with positive real parts and the two others eigenvalues are real).

\subsubsection{Eigenvalues at $I_{3}$($\rm{\textit{N}_3^*}$, $\rm{\textit{E}_3^*}$, $\rm{\textit{R}_3^*}$, 0)}

Concerning this last point $I_3$, the analytical analysis of stability is no more possible and it becomes then necessary to choose a parameter set.

\vspace{0.1in}

\textbf{Remark.} Let's notice that the number of real positive fixed points depends on the choice of parameters.

\vspace{0.1in}

\subsubsection{Bifurcation analysis}

Since $I_2$, $J_2$ and $I_3$ do not depend on parameters $(\beta_4, r_3, \alpha_1, \gamma_1)$, bifurcations can occur in these models (for both states of Colorado and Washington). Following the works of May \cite{May1976}, we propose to choose the parameter $\beta_1$, \textit{i.e.}, the growth rate of the population of the prey ($N$) or the rate of moving in and out of the Nonuser category, as bifurcation parameter. Let's notice that this choice is based on the assumption according to which the population of Nonusers increases in proportion of the demography. Then, for the same reasons as previously, an analytical analysis would be difficult even impossible. So, in the next section we will set all the parameters except $\beta_1$ and then we will use a bifurcation diagram to determine the values of the bifurcation parameters.

\subsubsection{Existence of bounded solutions}

Following the works of Freedman \cite{Freed} and while posing $N=x_1$, $E = x_2$, $R = x_3$ and $A = x_4$, the system (\ref{eq1}) can be rewritten as follows:

\begin{equation}
\label{eq8}
\begin{aligned}
\dfrac{dx_1}{dt} & = x_1 g \left( x_1 \right) - \left( r_1 x_2 + \alpha_1 x_4 + \alpha_2 x_3 \right) p_1\left( x_1 \right), \hfill \\
\dfrac{dx_2}{dt} & = x_2 \left( - \beta_2 + r_1 p_1\left( x_1 \right) \right) - \left(r_2 x_3 + \gamma_1 x_4 \right) p_2 \left( x_2 \right), \hfill \\
\dfrac{dx_3}{dt} & = x_3 \left(- \beta_3 + r_2 p_2 \left( x_2 \right) + \alpha_2 p_1 \left( x_1 \right) \right) - r_3 x_4 p_3 \left( x_3 \right) , \\
\dfrac{dx_4}{dt} & = x_4 \left[ - \beta_4 + r_3 p_3 \left( x_3 \right) + \alpha_1 p_1 \left( x_1 \right) + \gamma_1 p_2 \left( x_2 \right) \right],
\end{aligned}
\end{equation}

\vspace{0.1in}

where $g(x_1) = \beta_1 (1 - x_1)$, $p_i(x_i) = \dfrac{x_i}{h + x_i}$ with $i=1,2, 3$.

\vspace{0.1in}

Moreover, analysis of experimental data available on the prevalence of marijuana in the population of 21+ in the states of Colorado and Washington \cite{Hanley} has shown that the influence of $A$ on $N$ and $E$ as well as that of $R$ on $N$ are in fact very weak. So, we will consider that $\alpha_1 \ll 1$, $\alpha_2 \ll 1$ and $\gamma_1 \ll 1$. Thus, under these assumptions, model (\ref{eq8}) reads:

\begin{equation}
\label{eq9}
\begin{aligned}
\dfrac{dx_1}{dt} & = x_1 g \left( x_1 \right) - r_1 x_2 p_1\left( x_1 \right), \hfill \\
\dfrac{dx_2}{dt} & = x_2 \left[ - \beta_2 + r_1 p_1\left( x_1 \right) \right] - r_2 x_3 p_2 \left( x_2 \right), \hfill \\
\dfrac{dx_3}{dt} & = x_3 \left[- \beta_3 + r_2 p_2 \left( x_2 \right) \right] - r_3 x_4 p_3 \left( x_3 \right) , \\
\dfrac{dx_4}{dt} & = x_4 \left[ - \beta_4 + r_3 p_3 \left( x_3 \right) \right].
\end{aligned}
\end{equation}

Then, according to Theorem 2.1 stated by Freedman \cite[p. 72]{Freed}, ``all solutions initiating in the nonnegative cone are bounded and eventually enter a certain attracting set described below.''

\section{Applications of Predator-Prey NERA Model}

In order to perform numerical experiments for forecasting the marijuana consumption in the states of Colorado and Washington beyond the year of the I -- 502 implementation, we used data from Hanley \cite{Hanley}. The Washington State Institute for Public Policy conducted a benefit-cost analysis of the implementation of I -- 502, which legalizes recreational marijuana use for adults within the two states \cite{Hanley}. The data collected in the latter report is the result of an analysis of population-level data in order to monitor four key indicators of marijuana use, namely, current marijuana use, lifetime marijuana use, marijuana abuse or dependency and age of initiation prior to implementation of I -- 502. In the case of our NERA model, the first three categories correspond to the Recreational, Experimental and Addict category respectively \cite{Dauhoo}. The report highlights the importance of examining trends in that manner will allow them to monitor whether the implementation of I -- 502 appears to affect these key indicators of marijuana use over time. Our numerical experiments aim to forecast the marijuana consumption in the states of Colorado and Washington beyond the year of the I -- 502 implementation using data from Hanley \cite{Hanley}. The NERA predator-prey model is calibrated by estimating the parameters $r_1$, $r_2$, $r_3$, $\alpha_1$, $\alpha_2$, $\alpha_3$, $\beta_1$, $\beta_2$, $\beta_3$ and $\beta_4$. Parameter $h$ has been arbitrarily chosen equal to $1/2$ as usually done in theoretical ecology \cite{Scudo}. To this aim, we used these data and a genetic algorithm as explained in Dauhoo \textit{et al.} \cite{Dauhoo}. The fitness function used has been chosen to minimize the error that our model generates. MATLAB Optimtool is used and the fitness function is inserted in the genetic algorithm. We hence obtain the set of values in Tab. 2 and Tab. 3 corresponding to the NERA predator-prey model for Colorado and Washington. Obviously, since the parameters sets of both NERA predator-prey models have been calibrated starting from the data from Hanley \cite{Hanley}, numerical integration of Eqs. (\ref{eq1}) with parameters sets from Tab. 2 and 3 are in perfect agreement with the observed data. Thus, it confirms that the NERA predator-prey models can be a precious tool in forecasting of illicit drug consumption in the population of the state considered.

\begin{table}[h]
	\centering
	\caption{Parameter values for the consumption of Marijuana in Colorado}
\begin{tabular}{ccccccccccc}
  \hline
  $r_1$ & $r_2$ & $r_3$ & $\alpha_1$ & $\alpha_2$ & $\alpha_3$ & $\beta_1$ & $\beta_2$ & $\beta_3$ & $\beta_4$ \\
  \hline
  0.44 & 0.193 & 0.029 & 0.103 & 0.043 & 0.031 & 0.042 & 0.016 & 0.052 & 0.047 \\
  \hline
\end{tabular}
\label{tab2}
\end{table}	

\begin{table}[h]
	\centering
	\caption{Parameter values for the consumption of Marijuana in Washington}
\begin{tabular}{ccccccccccc}
  \hline
  $r_1$ & $r_2$ & $r_3$ & $\alpha_1$ & $\alpha_2$ & $\alpha_3$ & $\beta_1$ & $\beta_2$ & $\beta_3$ & $\beta_4$ \\
  \hline
  0.38 & 0.142 & 0.034 & 0.099 & 0.112 & 0.032 & 0.015 & 0.03 & 0.066 & 0.039 \\
  \hline
\end{tabular}
\label{tab3}
\end{table}	

\begin{figure}[htbp]
\centerline{\includegraphics[width = 12cm,height = 7.6cm]{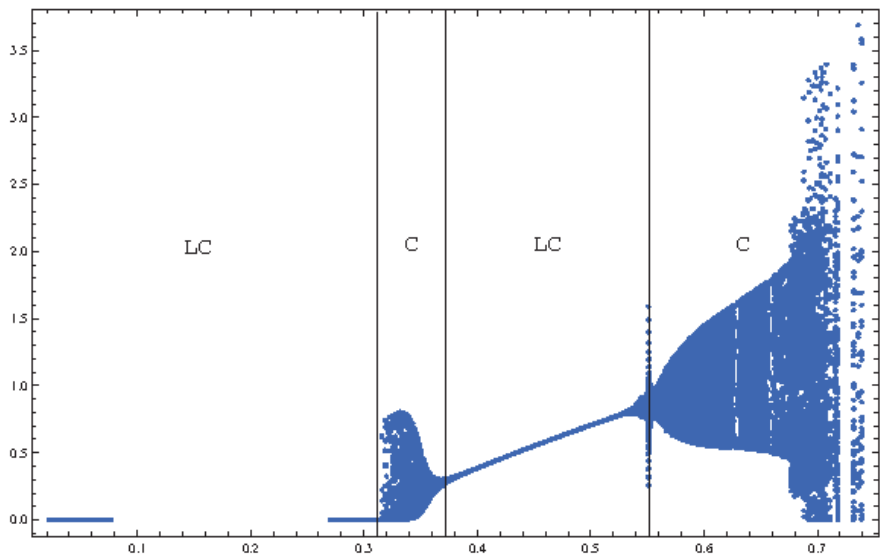}}
\caption{Bifurcation diagram of model (\ref{eq1}) $u_{max}$ as function of $\beta_1$.}
\label{fig2}
\end{figure}

\subsection{NERA model for Colorado}

Still using experimental data \cite{Hanley}, we set the parameters from Tab. \ref{tab2} where $\beta_1$ is the bifurcation parameter. By varying continuously $\beta_1$ from 0.02 to 0.8  (all other parameters are those given in Tab. 2), we determine the values for which bifurcations occur by plotting the bifurcation diagram presented in Fig. 2.

We observe in this diagram (Fig. 2), that for the value of $\beta_1 \in [0.02, 0.31]$ the solution of the NERA model (\ref{eq1}) for Colorado is a \textit{limit cycle} (LC) as exemplified in Fig. 3a. Then, still increasing parameter $\beta_1$ up to the first bifurcation value $\beta_1^{b_1} = 0.3175$, the solution becomes chaotic (C) as highlighted in Fig. 3b. Such chaotic feature persists up to the second bifurcation value $\beta_1^{b_2} = 0.37$ starting from which the solution becomes again a \textit{limit cycle} (LC) (see Fig 3c). Then, starting from the third bifurcation value $\beta_1^{b_3} = 0.55$ a chaotic attractor (C) appears again (see Fig. 3d).

\begin{figure}[htbp]
  \begin{center}
    \begin{tabular}{ccc}
      \includegraphics[width=6cm,height=6cm]{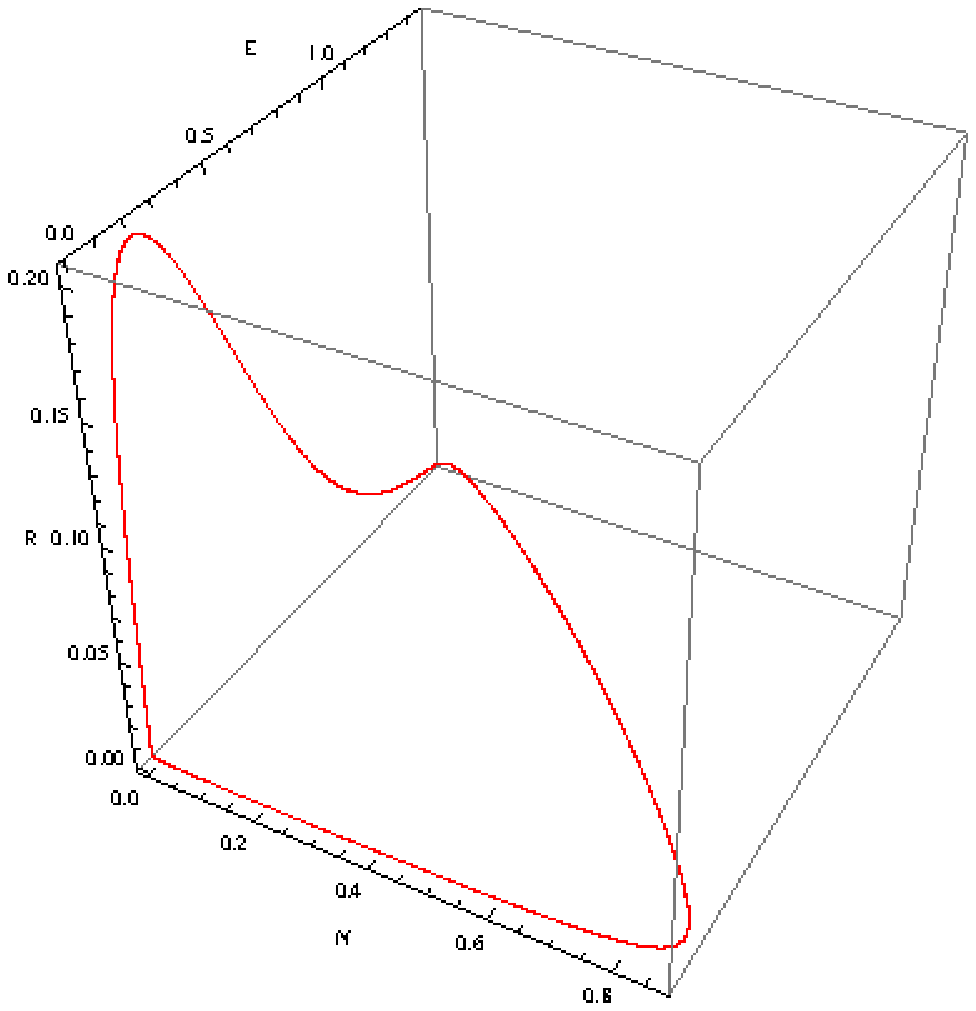} & ~~~~~~ &
      \includegraphics[width=6cm,height=6cm]{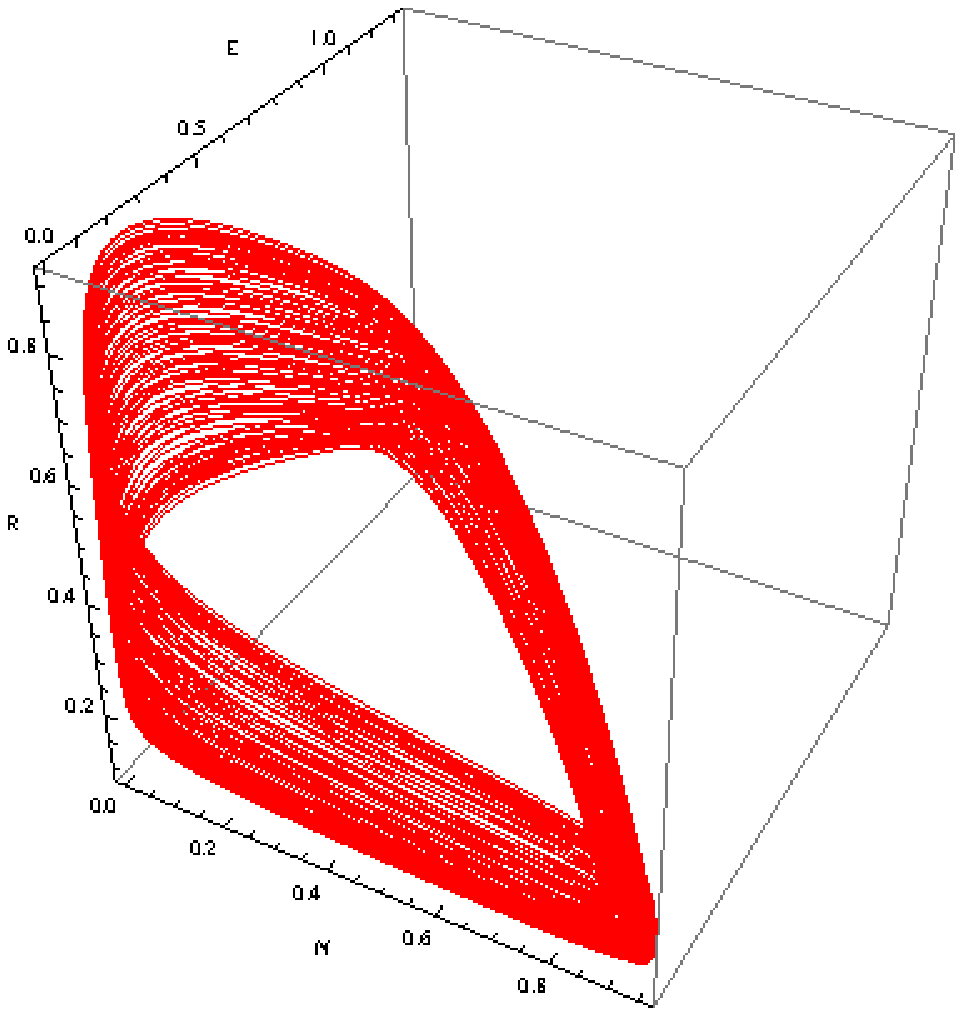} \vspace{0.1in} \\
      (a) $\beta_1=0.30$ & & (b) $\beta_1=0.35$ \\[0.2cm]
      \includegraphics[width=6cm,height=6cm]{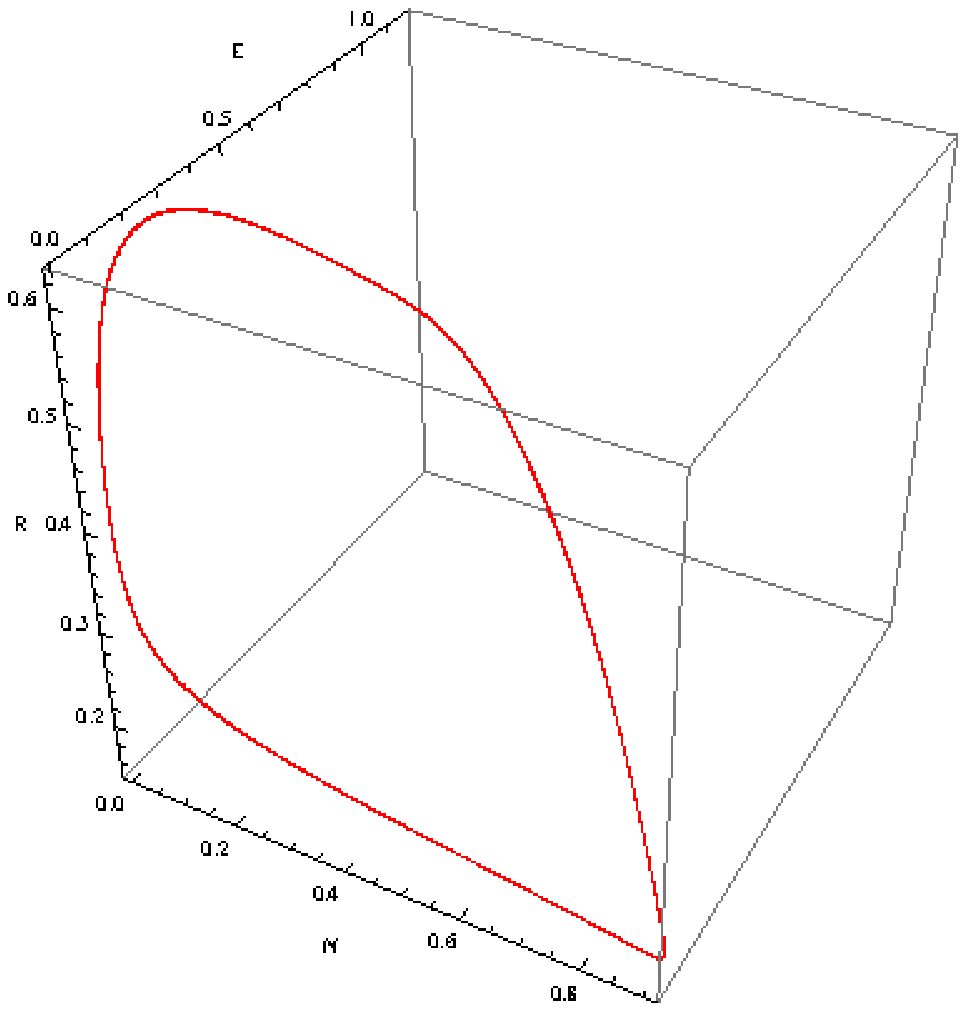} & ~~~ &
      \includegraphics[width=6cm,height=6cm]{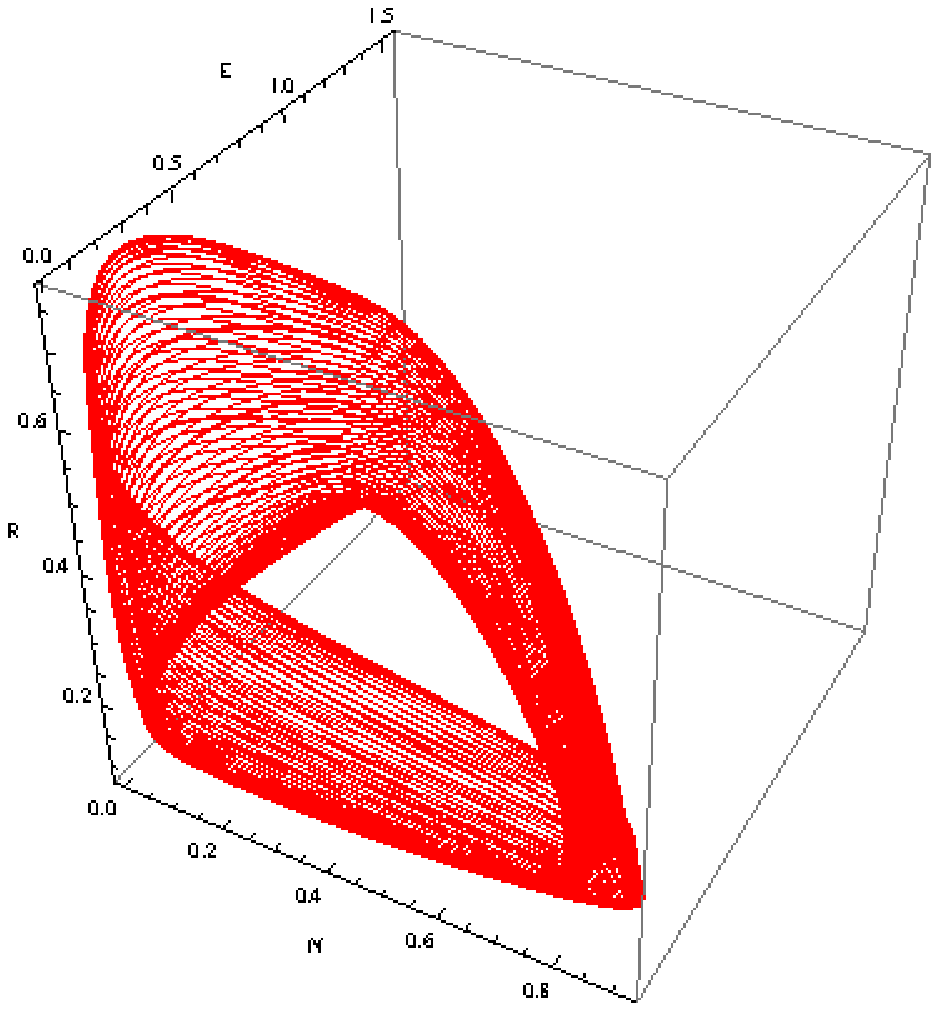} \vspace{0.1in} \\
      (c) $\beta_1 = 0.40$ & & (d) $\beta_1=0.60$ \\[-0.2cm]
    \end{tabular}
    \caption{Phase portraits of model (\ref{eq1}) in the ($N,E,R$)-space for various values $\beta_1$.}
    \label{fig3}
  \end{center}
  \vspace{-0.5cm}
\end{figure}

The algorithm developed by Marco Sandri \cite{Sandri} for Mathematica$^{\mbox{\scriptsize{\textregistered}}}$ has been used to perform the numerical calculation of the Lyapunov characteristics exponents (LCE) of the NERA predator-prey model (\ref{eq1}) for Colorado with the parameters from Tab. 2 and with $\beta_1  \in [0.02, 0.8]$. As an example for $\beta_1 = 0.3, 0.35, 0.4$ and $0.6$, Sandri's algorithm has provided respectively the following LCEs $(0, -0.0022, -0.013, -0.029)$, $(0, 0, -0.0010, -0.065)$, $(0, -0.0011, -0.0011, -0.066)$ and $(0, 0, -0.0017, -0.085)$. Then, according to the works of Klein and Baier \cite{KleinBaier}, a classification of (autonomous) continuous-time attractors of dynamical system (\ref{eq1}) on the basis of their Lyapunov spectrum, together with their Hausdorff dimension is presented in Tab. 4. LCEs values have been also computed with the Lyapunov Exponents Toolbox (LET) developed by Pr. Steve Siu for MatLab$^{\mbox{\scriptsize{\textregistered}}}$ and involving the two algorithms proposed by Wolf \textit{et al.} \cite{Wolf} and Eckmann and Ruelle \cite{Eckmann} (see https://fr.mathworks.com/matlabcentral/fileexchange/233-let). Results obtained by both algorithms are consistent.

\begin{table}[h]
    \centering
    \caption{LCEs of NERA model (\ref{eq1}) for Colorado for various values of $\beta_1$.}
\begin{tabular}{c c c c}\\[-2pt]
\hline
{\hspace{11mm} $\beta_1$} & LCE spectrum & Dynamics of the attractor & Hausdorff dimension  \\[6pt]
\hline\\[-2pt]
{\hspace{0.9mm} $0.02 < \beta_1 < 0.31$} & ($ 0, -, -, - $) & Periodic Motion (Limit Cycle) & $D = 1$ \\[1pt]
{\hspace{0.9mm} $0.31 < \beta_1 < 0.37$} & ($ 0, 0, -, - $) & 3-Torus (Quasi-Periodic Motion) & $D = 2$  \\[2pt]
{\hspace{0mm}   $0.37 < \beta_1 < 0.55$}  &  ($ 0, -, -, - $) & Periodic Motion (Limit Cycle) & $D = 1$ \\[1pt]
{\hspace{0.9mm} $0.55 < \beta_1 < 0.8$} & ($ 0, 0, -, - $) & 3-Torus (Quasi-Periodic Motion) & $D = 2$  \\[2pt]
\hline
\end{tabular}
\label{tab4}
\end{table}

\subsection{NERA model for Washington}

Now, we set the parameters from Tab. \ref{tab3} where $\beta_1$ is still the bifurcation parameter. By varying continuously $\beta_1$ from 0.02 to 0.36  (all other parameters are those given in Tab. 3), we determine the values for which bifurcations occur by plotting the bifurcation diagram presented in Fig. 4.

\begin{figure}[htbp]
\centerline{\includegraphics[width = 12cm,height = 7.6cm]{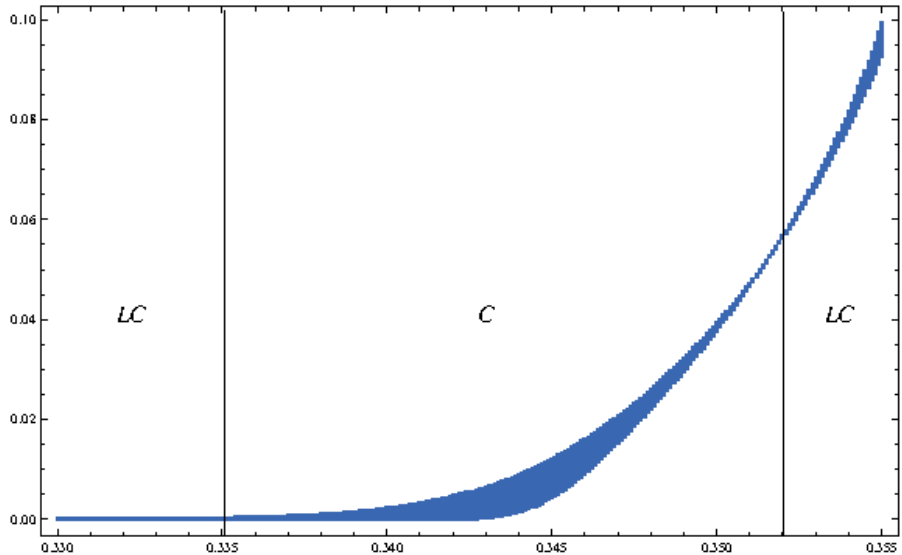}}
\caption{Bifurcation diagram of model (\ref{eq1}) $u_{max}$ as function of $\beta_1$.}
\label{fig4}
\end{figure}

We observe from this diagram (Fig. 4) that for the value of $\beta_1 \in [0.02, 0.334]$ the solution of the NERA model (\ref{eq1}) for Washington is a \textit{limit cycle} (LC) in the ($N,E$)-plane as exemplified in Fig. 5a. Then, still increasing parameter $\beta_1$ up to the first bifurcation value $\beta_1^{b_1} = 0.335$, the solution becomes chaotic (C) as highlighted in Fig. 5b \& Fig. 5c. Such chaotic feature persists up to the second bifurcation value $\beta_1^{b_2} = 0.357$ starting from which the solution becomes again a \textit{limit cycle} (see Fig 5d).

\begin{figure}[htbp]
  \begin{center}
    \begin{tabular}{ccc}
      \includegraphics[width=6cm,height=6cm]{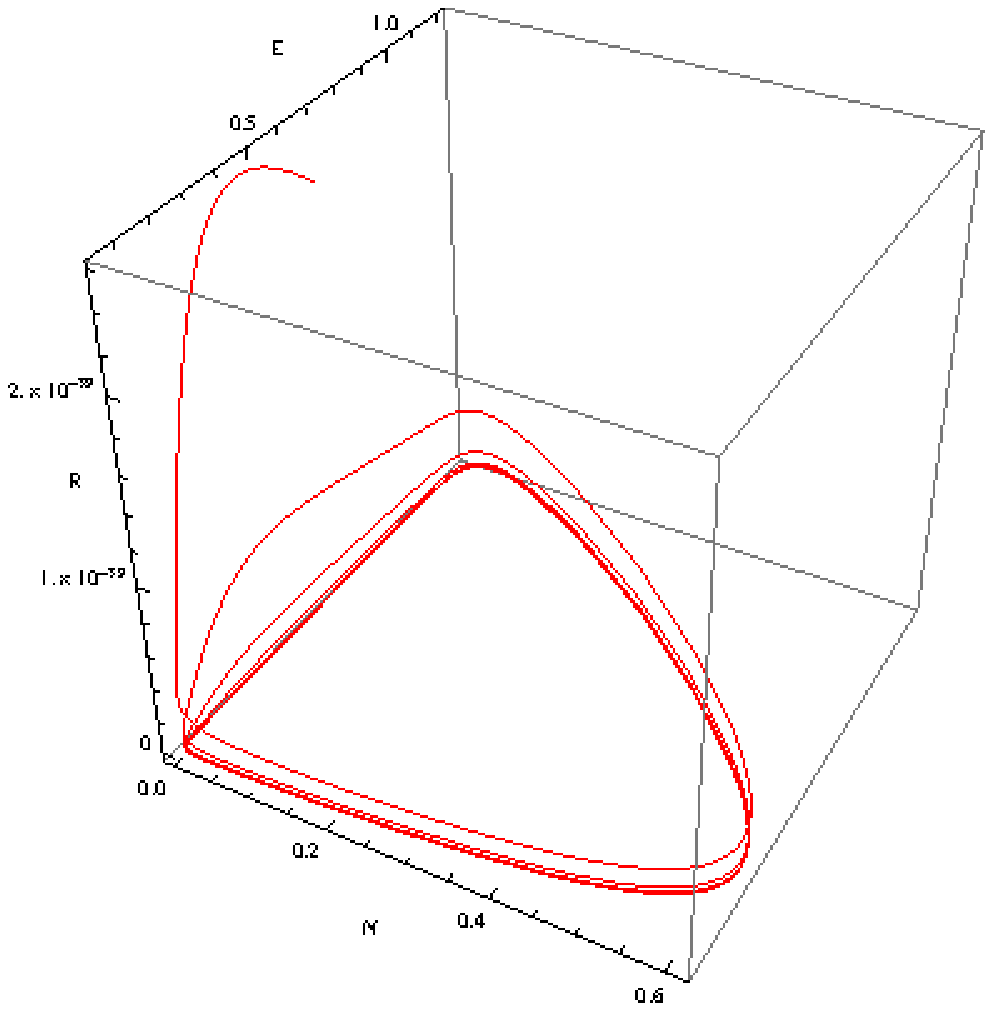} & ~~~~~~ &
      \includegraphics[width=6cm,height=6cm]{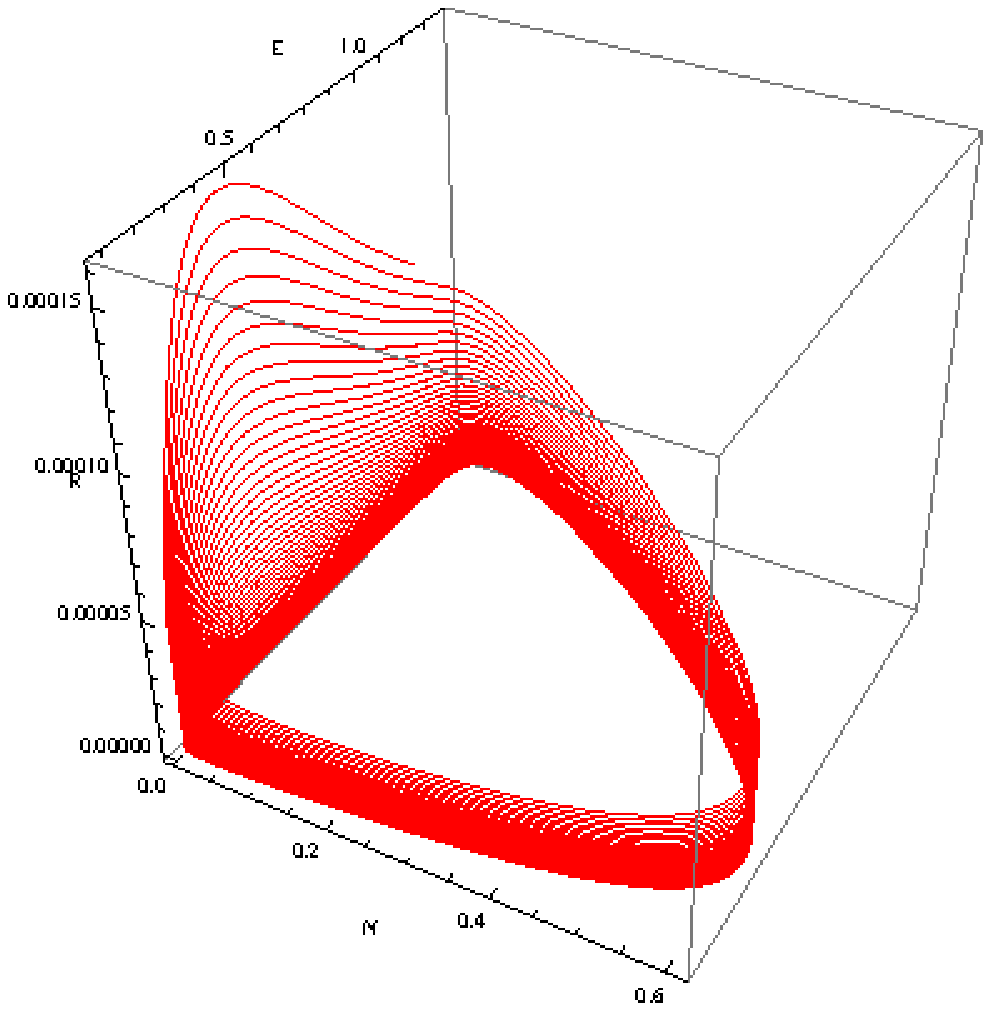} \vspace{0.1in} \\
      (a) $\beta_1=0.27$ & & (b) $\beta_1=0.34$ \\[0.2cm]
      \includegraphics[width=6cm,height=6cm]{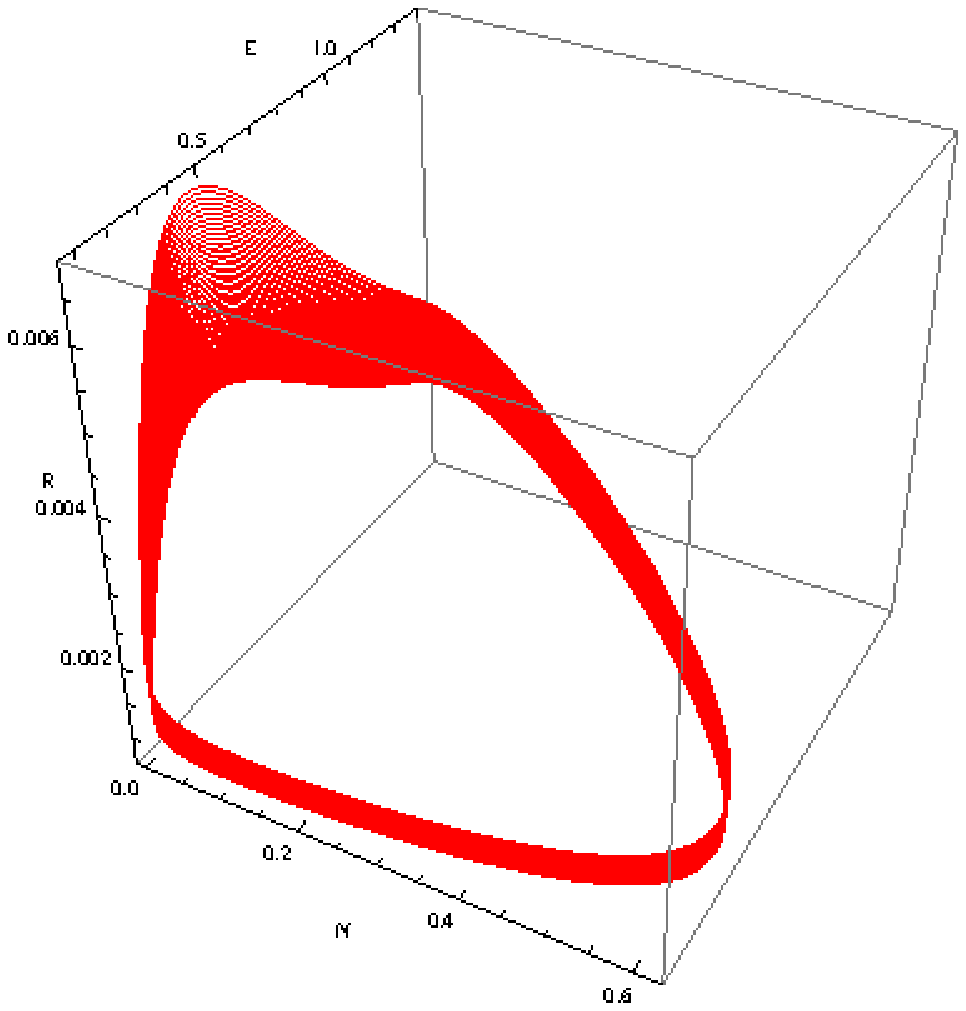} & ~~~ &
      \includegraphics[width=6cm,height=6cm]{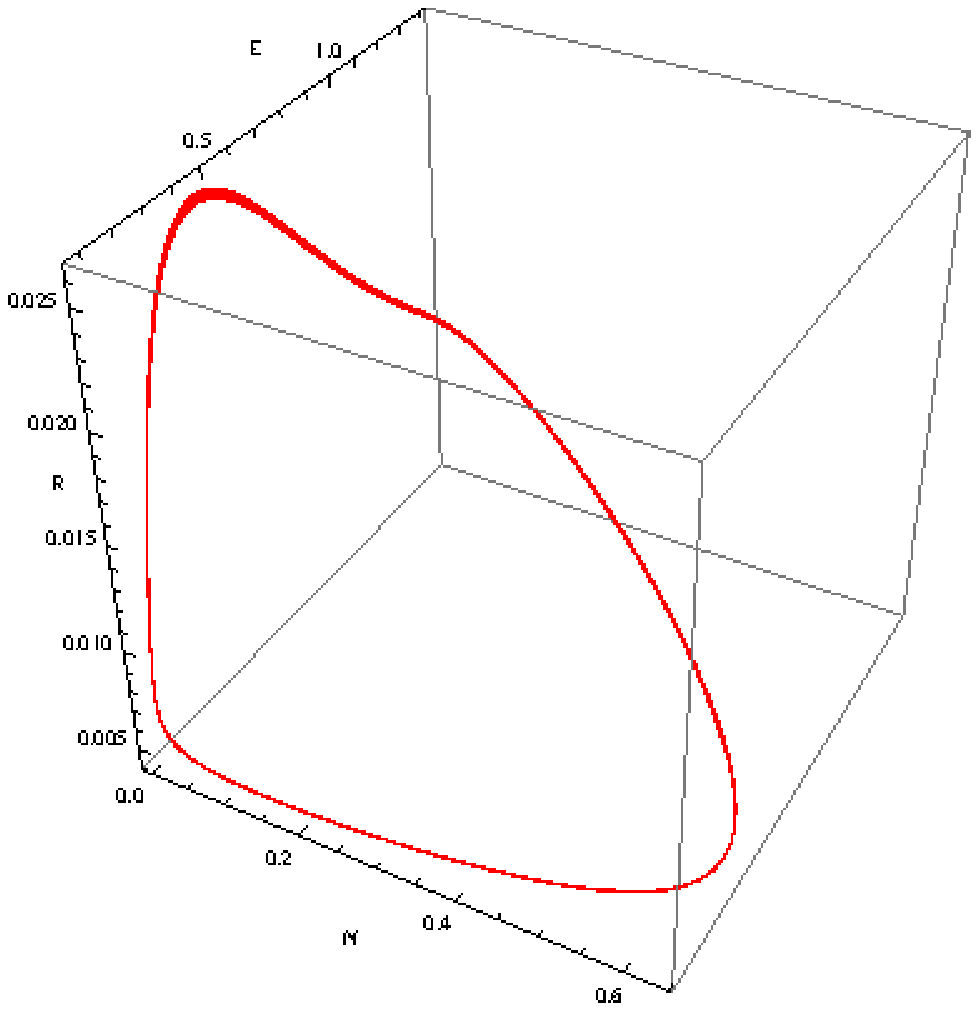} \vspace{0.1in} \\
      (c) $\beta_1 = 0.345$ & & (d) $\beta_1=0.3485$ \\[-0.2cm]
    \end{tabular}
    \caption{Phase portraits of model (\ref{eq1}) in the ($N,E,R$)-space for various values $\beta_1$.}
    \label{fig5}
  \end{center}
  \vspace{-0.5cm}
\end{figure}

Still using the algorithm developed by Marco Sandri \cite{Sandri} we numerically compute the Lyapunov characteristics exponents (LCE) of the NERA model (\ref{eq1}) for Washington with the parameters from Tab. 3 and with $\beta_1  \in [0.02, 0.36]$. In this case, for $\beta_1 = 0.27, 0.34, 0.345$ and $0.3485$, Sandri's algorithm provides respectively the following LCEs $(0, -0.0078, -0.0089, -0.019)$, $(0, 0, -0.010, -0.017)$, $(0, 0, -0.0092, -0.016)$ and $(0, 0, -0.011, -0.016)$.\\ Then, as previously, a classification of attractors of dynamical system (\ref{eq1}) on the basis of its Lyapunov spectrum, together with its Hausdorff dimension is presented in Tab. 5. Here again, LCEs values have been also computed with the Lyapunov Exponents Toolbox (LET) and results obtained by both algorithms are consistent.

\begin{table}[h]
    \centering
    \caption{LCEs of NERA model (\ref{eq1}) for Washington for various values of $\beta_1$.}
\begin{tabular}{c c c c}\\[-2pt]
\hline
{\hspace{11mm} $\beta_1$} & LCE spectrum & Dynamics of the attractor & Hausdorff dimension  \\[6pt]
\hline\\[-2pt]
{\hspace{0.9mm} $0.2 < \beta_1 < 0.33$} & ($ 0, -, -, - $) & Periodic Motion (Limit Cycle) & $D = 1$ \\[1pt]
{\hspace{0.9mm} $0.3375 < \beta_1 < 0.44$} & ($ 0, 0, -, - $) & 3-Torus & $D = 2$  \\[2pt]
{\hspace{0mm}   $0.44 < \beta_1 < 0.45$}  &  ($ 0, -, -, - $) & Periodic Motion (Limit Cycle) & $D = 1$ \\[1pt]
\hline
\end{tabular}
\label{tab5}
\end{table}

\newpage

\section{Discussion}

By considering that drug users are classified into four main categories: non (N), experimental (E), recreational (R) and addicts (A) users, Dauhoo \textit{et al.} \cite{Dauhoo} have proposed the NERA model. Nevertheless, although this four-dimensional dynamical system took into account the mutual influence that drug users (E, R and A) can have on non-users (N) and on each other, it did not contain any limitation in the growth and decay of each variable. Thus, no oscillatory or chaotic regime could be observed. So, the aim of this work was to propose a modified version of this NERA model by analogy with the classical predator-prey models and while considering non-users (N) as prey and users (E, R and A) as predator. Thus, this new model included a ``stabilizing effect'' of the growth rate of the preys (N) and a ``destabilizing effect'' of the predators (E, R and A) saturation. Functional responses of Verhulst and Holling type II have been used for modeling these effects. Then, in order to perform numerical experiments for forecasting the marijuana consumption in the states of Colorado and Washington beyond the year of the I -- 502 implementation, we used data from Hanley \cite{Hanley} and a genetic algorithm as explained in Dauhoo \textit{et al.} \cite{Dauhoo}. Thus, the NERA predator-prey model has been calibrated by estimating all the parameters except $h$ which has been arbitrarily chosen equal to $1/2$ as usually done in theoretical ecology \cite{Scudo}. We hence obtained two parameters sets corresponding to the NERA predator-prey model for the state of Colorado and Washington. Following the works of May \cite{May1976}, we chose the parameter $\beta_1$, \textit{i.e.}, the growth rate of the population of the prey ($N$) or the rate of moving in and out of the Nonuser category, as bifurcation parameter. This choice was based on the assumption that the population of Nonusers increases in proportion of the demography. A stability and bifurcation analysis of the NERA model for Colorado and for Washington has therefore been performed.

Concerning the NERA model for Colorado, the bifurcation diagram has shown that when the value of parameter $\beta_1 \in [0.02, 0.31]$, the solution is a \textit{limit cycle} confirming thus the behavior of the observed data from Hanley \cite{Hanley}. So, the number of individuals of each group N, E, R and A oscillates in a deterministic way with a period and amplitude that can be numerically computed. Then, still increasing parameter $\beta_1$ up to the first bifurcation value $\beta_1^{b_1} = 0.3175$, we have shown that the solution becomes quasi periodic (3-Torus). Such chaotic attractor persists up to the second bifurcation value $\beta_1^{b_2} = 0.37$ starting from which the solution becomes again a \textit{limit cycle}. When parameter $\beta_1^{b_3} = 0.55$ reaches the third bifurcation value a chaotic attractor appears again. These results have been confirmed by the computation of the Lyapunov characteristics exponents.

Concerning the NERA model for Washington, the bifurcation diagram has shown that when the value of parameter $\beta_1 \in [0.02, 0.26]$ the solution is a \textit{limit cycle} confirming thus the behavior of the observed data from Hanley \cite{Hanley}. Then, still increasing parameter $\beta_1$ up to the first bifurcation value $\beta_1^{b_1} = 0.34$, the solution becomes quasi periodic (3-Torus). Such chaotic attractor persists up to the second bifurcation value $\beta_1^{b_2} = 0.357$ starting from which the solution becomes again a \textit{limit cycle}. These results have been also confirmed by the computation of the Lyapunov characteristics exponents.

Thus, we have shown that an increase of the population of non-users (N) leads first to a periodic evolution of drug-users (E, R, A). But when this increase reaches a certain threshold for both states of Colorado and Washington, the evolution of all variables becomes unpredictable. Such result can be also interpreted by analogy with the so-called ``paradox of enrichment'' which states that increasing the food available to the prey caused the predator's population to destabilize \cite{Rosenzweig}. We have thus confirmed on one hand that the NERA predator-prey models can be a precious tool in forecasting of illicit drug consumption in the population of the state considered and, on the other hand that they can be of substantial interest to policy-makers in the problematic of illicit drug consumption.

\acknowledgments
Matja{\v z} Perc was supported by the Slovenian Research Agency (Grants J1-7009, J4-9302, J1-9112 and P5-0027)

\end{document}